\documentclass[a4paper]{jpconf}
\usepackage{amsfonts}
\usepackage{graphicx}

\newtheorem{theorem}{Theorem}
\newtheorem{prop}[theorem]{Proposition}
\newtheorem{lem}[theorem]{Lemma}
\newtheorem{unnumber}{}

\newtheorem{prepf}[unnumber]{Proof}
\newenvironment{pf}{\prepf\rm}{\endprepf}
\newcommand{\qed}{$\Box$}

\newcommand{\ee}{\mathord{\mathrm{e}}}

\def\proof{\par\noindent{\bf Proof.\enspace}\rm}
\def\ox{\overline{x}}
\def\oy{\overline{y}}
\def\oR{\overline{R}}
\def\bp{B^\prime}
\def\BP{{\mathbb P}}
\def\BE{{\mathbb E}}

\begin{document}

\title{Asymptotic enumeration of incidence matrices}

\author{Peter Cameron, Thomas Prellberg and Dudley Stark}

\address{School of Mathematical Sciences, Queen Mary, University of London, Mile End Road, London E1 4NS, United Kingdom}

\ead{p.j.cameron, t.prellberg, d.stark@qmul.ac.uk}

\begin{abstract}
We discuss the problem of counting {\em incidence matrices},
i.e. zero-one matrices with no zero rows or columns.
Using different approaches we give three different proofs for
the leading asymptotics for the number of matrices with 
$n$ ones as $n\to\infty$. We also give refined results for the
asymptotic number of $i\times j$ incidence matrices with $n$ ones.
\end{abstract}.

\section{Introduction}

We call an {\em incidence matrix} a zero-one matrix with no zero rows and
columns and denote by $F(n)$ the number of incidence matrices with exactly 
$n$ ones, where $n\in\mathbb N$. For example, the four 
incidence matrices with $n=2$ are \[\pmatrix{1&1\cr},\quad \pmatrix{1\cr1\cr},\quad
\pmatrix{1&0\cr0&1\cr},\quad \pmatrix{0&1\cr1&0\cr}.\]
The first few terms of the sequence $F(n)$ for $n\in\mathbb N$ are
$$1,4,24,196,2016,24976,361792,5997872,111969552,\ldots$$
taken from the \emph{On-Line Encyclopedia of Integer
Sequences}~\cite{oeis}, where this appears as sequence A101370. For convenience,
we further define $F(0)=1$.

If one imposes additional symmetries or constraints, such as allowing or
prohibiting repeated rows or columns, or considering equivalence classes
under row or column permutations, one is led to many different enumeration
problems, as discussed in \cite{cps}.

The counting problem can be interpreted in a surprisingly rich variety of
different ways, leading to rather different mathematical approaches.

\begin{itemize}

\item{\bf Counting hypergraphs by weight}

Given a hypergraph on the vertex set $\{x_1,\ldots,x_r\}$, with edges
$E_1,\ldots,E_s$ (each a non-empty set of vertices), the \emph{incidence
matrix} $A=(a_{ij})$ is the matrix with $(i,j)$ entry $1$ if $x_i\in E_j$,
and $0$ otherwise. The \emph{weight} of the hypergraph is the sum of the
cardinalities of the edges. Thus $F(n)$ is the number of vertex- and edge-labelled
hypergraphs of weight~$n$ with no isolated vertices, up to isomorphism.

\item{\bf Counting bipartite graphs by edges}

Given a zero-one matrix $A=(A_{ij})$, there is a (simple) bipartite graph
whose vertices are indexed by the rows and columns of $A$, with an edge
from $r_i$ to $c_j$ if $A_{ij}=1$. The graph has a distinguished bipartite
block (consisting of the rows). Thus, $F(n)$ counts
labelled bipartite graphs with $n$ edges and a distinguished
bipartite block.

\item{\bf Counting pairs of partitions, or binary block designs}

A \emph{block design} is a set of \emph{plots} carrying two partitions,
the \emph{treatment partition} and the \emph{block partition}. It is said
to be \emph{binary} if no two distinct points lie in the same part of both
partitions; that is, if the meet of the two partitions is the partition into
singletons. Thus, $F(n)$ is the number of binary block designs with
$n$ plots and labelled treatments and blocks.

\item{\bf Counting orbits of certain permutation groups}

A permutation group $G$ on a set $X$ is \emph{oligomorphic} if the number
$F^*_n(G)$ of orbits of $G$ on $n$-tuples of elements of $X$ is finite for all 
$n$. Equivalently, the number $F_n(G)$ of orbits on ordered $n$-tuples of 
distinct elements of $X$ is finite, and the number $f_n(G)$ of orbits on 
$n$-element subsets of $X$ is finite, for all $n$. These numbers satisfy 
various conditions, including the following:
  \begin{itemize}
  \item[$\bullet$] $F^*_n(G)=\sum_{k=1}^nS(n,k)F_k(G)$ and its
inverse $F_n(G)=\sum_{k=1}^ns(n,k)F_k^*(G)$, with $s(n,k)$ and 
$S(n,k)$ Stirling numbers of the first and second kind, respectively;
  \item[$\bullet$] $f_n(G)\le F_n(G)\le n!f_n(G)$, where the right-hand bound is attained
if and only if the group induced on a finite set by its setwise stabiliser
is trivial.
  \end{itemize}

For example, let 
$A$ be the group of all order-preserving permutations of the rational numbers. 
Then
$f_n(A)=1$ and $F_n(A)=n!$~.

Now if $H$ and $K$ are permutation groups on sets $X$ and $Y$, then the direct
product $H\times K$ acts coordinatewise on the Cartesian product $X\times Y$.
It is easy to see that $F_n^*(H\times K)=F_n^*(H)F_n^*(K)$.

Let $(x_1,y_1)$, \dots, $(x_n,y_n)$ be $n$ distinct elements of $X\times Y$.
If both $X$ and $Y$ are ordered, then the set of $n$ pairs can be described by
a matrix with $n$ ones in these positions, where the rows and columns of the
matrix are indexed by the sets $\{x_1,\ldots,x_n\}$ and $\{y_1,\ldots,y_n\}$
respectively (in the appropriate order).
Thus \[F(n)=f_n(A\times A).\]
Discussion of this ``product action'' can be found in~\cite{cgm}.

\end{itemize}

For an extended discussion of these interpretations see \cite{cps}.
For instance, when considering hypergraphs it is more natural to consider the unlabelled problem,
which leads to identification of incidence matrices which are equivalent under
permutation of rows or columns. Also, forbidding repeated rows corresponds to
counting simple hypergraphs with no repeated edges.

\section{The asymptotics of $F(n)$}

It is possible to compute $F(n)$ explicitly. For fixed~$n$, let
$m_{ij}(n)$ be the number of $i\times j$ matrices with $n$ ones (and no zero
rows or columns). We set $m_{00}(0)=1$ and $F(0)=1$. Then
\begin{equation}\label{firstdisp}
\sum_{i\le k}\sum_{j\le l}{k\choose i}{l\choose j}m_{ij}(n)={kl\choose n},
\end{equation}
so by M\"obius inversion,
\begin{equation}\label{mobius}
m_{kl}(n)=\sum_{i\le k}\sum_{j\le l}(-1)^{k+l-i-j}{k\choose i}{l\choose j}
{ij\choose n},
\end{equation}
and then
\begin{equation}\label{fsum}
F(n)=\sum_{i\le n}\sum_{j\le n}m_{ij}(n).
\end{equation}

For sequence $a_n$, $b_n$, we use the notation $a_n\sim b_n$
to mean $\lim_{n\to\infty}a_n/b_n=1$.
It is clear from the argument above that
\[F(n)\le{n^2\choose n}\sim\frac{1}{\sqrt{2\pi n}}(n\ee)^n,\]
and of course considering permutation matrices shows that
\[F(n)\ge n!\sim\sqrt{2\pi n}\left(\frac{n}{\ee}\right)^n.\]
\begin{theorem}\label{F1111}
$$
F(n)\sim
\frac{n!}4e^{-\frac12(\log2)^2}\frac1{(\log2)^{2n+2}}\;.
$$
\end{theorem}

We remark that for $n=10$, the asymptotic expression is about $2.5\%$
less than the actual value of $2324081728$.

As announced in \cite{cps}, we have three different proofs of 
Theorem \ref{F1111}. The first proof employs pairs of random preorders
and a probabilistic argument, the second proof uses counting of orbits 
of products of permutation groups, and the third 
proof employs a surprisingly simple identity. 

\noindent {\bf First Proof:}
This proof uses a procedure which, when successful,
generates an incidence matrix uniformly
at random from all incidence matrices. The probability of success can
be estimated and the asymptotic formula for $F(n)$ results.

Let $R$ be a binary relation on a set $X$. We say $R$ is {\em reflexive}
if $(x,x)\in R$ for all $x\in X$. We say $R$ is {\em transitive}
if $(x,y)\in R$ and $(y,z)\in R$ implies $(x,z)\in R$.
A {\em partial preorder} is a relation $R$ on $X$ which is reflexive
and transitive. A relation $R$ is said to satisfy {\em trichotomy}
if, for any $x,y\in X$, one of the cases $(x,y)\in R$, $x=y$, or $(y,x)\in R$
holds. We say that $R$ is a {\em preorder} if it is a partial preorder
that satisfies trichotomy. The members of $X$ are said to be the
{\em elements} of the preorder.

A relation $R$ is {\em antisymmetric} if, whenever $(x,y)\in R$ and
$(y,x)\in R$ both hold, then $x=y$. A relation $R$ on $X$ is a
{\em partial order} if it is reflexive, transitive, and antisymmetric.
A relation is a {\em total order}, if it is a partial order
which satisfies trichotomy.
Given a partial preorder $R$ on $X$,
define a new relation $S$ on $X$ by the rule that
$(x,y)\in S$ if and only if both $(x,y)$ and $(y,x)$ belong to $R$. Then
$S$ is an equivalence relation. Moreover, $R$ induces a partial order
$\ox$ on the set of equivalence classes of $S$ in a natural way:
if $(x,y)\in R$, then $(\ox,\oy)\in\oR$, where $\ox$ is the
$S$-equivalence class containing $x$ and similarly for $y$.
We will call an $S$-equivalence class a {\em block}.
If $R$ is a preorder,
then the relation $\oR$ on the equivalence classes of $S$ is a total order.
See Section 3.8 and question 19 of Section 3.13 in \cite{C}
for more on the above definitions and results.
Random preorders are considered in \cite{cs}.

Given a preorder on elements $[n]:=\{1,2,\ldots,n\}$
with $K$ blocks, let $B_1,B_2,\ldots,B_K$
denote the blocks of the preorder. Generate two preorders uniformly
at random, $B_1,B_2,\ldots,B_K$ and $\bp_1,\bp_2,\ldots,\bp_L$.
For each $1\leq i<j\leq n$, define the event $D_{i,j}$ to be
$$
D_{i,j}=\{{\rm for \ each \ of \ the \ two \ preorders \ }
i {\rm \ and \ } j {\rm \ are \ in \ the \ same \ block}\}.
$$
Furthermore, define
$$
W=\sum_{1\leq i<j\leq n}I_{D_{i,j}},
$$
where the indicator random variables are defined by
$$
I_{D_{i,j}}=
\left\{
\begin{array}{l l}
1&{\rm \  if \ } D_{i,j} {\rm \ occurs};\\
0&{\rm otherwise}.
\end{array}
\right.
$$
If $W=0$, then the procedure is successful, in which case
$B_k\cap\bp_l$ consists of either 0 or 1 elements for each
$1\leq k\leq K$ and $1\leq l\leq L$. If the procedure is successful,
then we define
the corresponding $K\times L$ incidence matrix $A$ by
$$
A_{k,l}=
\left\{
\begin{array}{l l}
1&{\rm \  if \ } B_k\cap \bp_l\neq\emptyset;\\
0&{\rm \  if \ } B_k\cap \bp_l=\emptyset.
\end{array}
\right.
$$
It is easy to check that the above definition of $A$ in fact produces
an incidence matrix and that each incidence matrix occurs in $n!$ different 
ways by the construction. It follows that
$$
F(n)=\frac{P(n)^2\BP(W=0)}{n!},
$$
where $P(n)$ is the number of preorders on $n$ elements if $n\geq 1$
and $P(0)=1$.

It is known (see \cite{Bar}, for example)
that the exponential generating function of $P(n)$
is
\begin{equation}\label{pgenfunc}
\sum_{n=0}^\infty\frac{P(n)}{n!}z^n=\frac{1}{2-e^z}.
\end{equation}
The preceding equality implies that
$P(n)$ has asymptotics given by
\begin{equation}\label{Pasymp}
P(n)\sim \frac{n!}{2}\left(\frac{1}{\log2}\right)^{n+1}\;.
\end{equation}
It remains to find the asymptotics of $\BP(W=0)$.

The $r$th falling moment of $W$ is
\begin{eqnarray}
\BE(W)_r&=& \BE W(W-1)\cdots(W-r+1)\nonumber\nonumber\\
&=&
\BE\left(\sum_{{\rm pairs \ }(i_s,j_s){\rm \ different}}
I_{i_1,j_1}\cdots I_{i_r,j_r}\right)\label{falling}\\
&=&
\BE\left(\sum_{{\rm all \ }i_s {\rm \ and \ } j_s{\rm \ different}}
I_{i_1,j_1}\cdots I_{i_r,j_r}\right)+
\BE\left(\sum^\ast
I_{i_1,j_1}\cdots I_{i_r,j_r}\right),\label{falling2}
\end{eqnarray}
with ${\displaystyle \sum^\ast}$ defined
to be the sum with all pairs $(i_s,j_s)$ different, but not all
$i_s, j_s$ different.

First we find the asymptotics of the first term in (\ref{falling2}).
For given sequences $i_1,i_2,\ldots,i_r$, $j_1,j_2,\ldots,j_r$,
the expectation
$\BE(I_{i_1,j_1}\cdots I_{i_r,j_r})$ is
the number of ways
of forming two preorders on the set of elements
$[n]\setminus\{j_1,j_2,\ldots,j_r\}$ and then for each $s$
adding the element
$j_s$ to the block containing $i_s$ in both preorders
(which ensures that $D_{i_s,j_s}$ occurs for each $s$) and dividing the
result by $P(n)^2$.
Since the number of ways of choosing  
$i_1,i_2,\ldots,i_r$, $j_1,j_2,\ldots,j_r$ equals
$\frac{n!}{2^r(n-2r)!}$,
This gives
\begin{eqnarray*}
\BE\left(\sum_{{\rm all \ }i_s {\rm \ and \ } j_s{\rm \ different}}
I_{i_1,j_1}\cdots I_{i_r,j_r}\right)&=&
\frac{n!}{2^r(n-2r)!}
\frac{P(n-r)^2}{P(n)^2}\\
&\sim&
\left(\frac{(\log 2)^2}{2}\right)^r,
\end{eqnarray*}
where we have used (\ref{Pasymp}).

The second term is bounded in the following way. For each sequence
$(i_1,j_1), (i_2,j_2),\ldots, (i_s,j_s)$ in the second term we form
the graph $G$ on vertices $\bigcup_{s=1}^r \{i_s,j_s\}$ with edges
$\bigcup_{s=1}^r \{\{i_s,j_s\}\}$.
Consider the unlabelled graph $G^\prime$ corresponding to $G$
consisting of $v$ vertices and $c$ components.
The number of ways
of labelling $G^\prime$ to form $G$ is bounded by
$n^v$. The number of preorders corresponding
to this labelling is $P(n-v+c)$ because we form a preorder on $n-v+c$ vertices
after which
the vertices in the connected component of $G$ containing a particular
vertex get added to that block.
Therefore, we have
\begin{eqnarray*}
\BE\left(\sum^\ast
I_{i_1,j_1}\cdots I_{i_r,j_r}\right)&\leq&
\sum_{G^\prime} n^v \frac{P(n-v+c)^2}{P(n)^2}\\
&=&\sum_{G^\prime}O\left(n^{2c-v}\right),
\end{eqnarray*}
where the constant in $O\left(n^{2c-v}\right)$ is uniform over all $G^\prime$
because $v\leq 2r$.
Since at least one vertex is adjacent
to more than one edge, the graph $G$ is not a perfect matching.
Furthermore, each component of $G$ contains at least two vertices.
It follows that $2c<v$ and,
as a result,
$$
\BE\left(\sum^\ast
I_{i_1,j_1}\cdots I_{i_r,j_r}\right)=O\left(n^{-1}\right).
$$

The preceding analysis shows that
$$
\BE(W)_r\sim\left(\frac{(\log 2)^2}{2}\right)^r
$$
for each $r\geq 0$. The method of moments implies that the distribution
converges weakly to the distribution of a Poisson$((\log 2)^2/2)$ distributed
random variable and therefore
\begin{equation}\label{W0}
\BP(W=0)\sim\exp\left(-\frac{(\log 2)^2}{2}\right).
\end{equation}
\qed

\noindent {\bf Second Proof:}
We now give a proof using product actions of groups, as discussed in the 
introduction.  First of all, this approach leads to a 
different and simpler expression than (\ref{fsum}) 
for $F(n)$ as a sum of terms of alternating sign.

\begin{prop}\label{prop2}
\[F(n)=\frac{1}{n!}\sum_{k=1}^ns(n,k)P(k)^2,\]
where
\[P(n)=\sum_{k=1}^nS(n,k)k!\]
is the number of (total) preorders of $\{1,\ldots,n\}$, and $s(n,k)$ and
$S(n,k)$ are Stirling numbers of the first and second kind respectively.
\end{prop}

This is proved in~\cite{cgm}, but can be seen as follows. Using
the group $A$ of all order-preserving permutation groups acting on $\mathbb Q$,
we consider the direct product $A\times A$ acting on $\mathbb Q\times\mathbb Q$.
We have $F_n(A)=n!$, whence it follows that $F^*_n(A)=\sum_{k=1}^nS(n,k)k!=P(n)$.
Thus
\[F^*_n(A\times A)=P(n)^2=\sum_{k=1}^n S(n,k)F_k(A\times A),\]
and so the inverse relation between the two kinds of Stirling numbers gives
\[F_n(A\times A)=\sum_{k=1}^ns(n,k)P(k)^2.\]
Finally, the group $A\times A$ has the property that the setwise stabiliser
of a finite set fixes it pointwise, and so $f_n(A\times A)=F_n(A\times A)/n!$~.\qed

We now replace $P(k)$ by the asymptotic form (\ref{Pasymp}) given earlier.
For $k\ge n/2$, the difference is
exponentially small; and we will show below that the contribution of the
terms with $k<n/2$ is negligible, so it suffices to note that the error
we make is smaller than the approximated term.

So let
\[F'(n) = \frac{1}{4}\cdot\frac{1}{n!}\sum_{k=1}^n s(n,k)(k!)^2c^{k+1},\]
where $c=1/(\log2)^2$ is as in the statement of the theorem. As we have argued,
$F(n)\sim F'(n)$.

Now $(-1)^{n-k}s(n,k)$ is the number of permutations in the symmetric group
$S_n$ which have $k$ cycles. So we can write the formula for $F'(n)$ as a sum
over~$S_n$, where the term corresponding to a permutation with $k$ cycles is
$(-1)^{n-k}(k!)^2c^{k+1}$. In particular, the identity permutation gives
us a contribution
\[g(n) = \frac{1}{4}\,n!\,c^{n+1},\]
and we have to show that $F'(n)\sim Cg(n)$ as $n\to\infty$, where
$C=\exp(-(\log2)^2/2)$.

To prove this, we write $F'(n)=F'_1(n)+F'_2(n)+F'_3(n)$, where the three
terms are sums over the following permutations:
\begin{description}
\item{$F'_1$:} all involutions (permutations with $\sigma^2=1$);
\item{$F'_2$:} the remaining permutations with $k\ge\lceil n/2\rceil$;
\item{$F'_3$:} the rest of $S_n$.
\end{description}
We argue that $F'_1(n)\sim Cg(n)$, while $F'_2(n),F'_3(n)=o(g(n))$.

\subparagraph{Case $F'_1$:} Let $l=n-k$. Now an involution with $k$ cycles
has $l$ cycles of length~$2$ and $n-2l$ fixed points; so $l\le n/2$. The
number of such permutations is
\[{n\choose 2l}\frac{(2l)!}{2^l\,l!}=\frac{n(n-1)\cdots(n-2l+1)}{2^l\,l!}.\]
So
\begin{eqnarray*}
\frac{F'_1(n)}{g(n)}
&=& \sum_{l=0}^{\lfloor n/2\rfloor}\frac{n(n-1)\cdots(n-2l+1)}{2^l\,l!}(-1)^l
\frac{((n-l)!)^2}{(n!)^2}c^{-l}\\
&=& \sum_{l=0}^{\lfloor n/2\rfloor}\frac{1}{l!}\left(\frac{-1}{2c}\right)^l
\frac{(n-l)\cdots(n-2l+1)}{n\cdots(n-l+1)}.
\end{eqnarray*}

Now
\[\sum_{l=0}^\infty \frac{1}{l!}\left(\frac{-1}{2c}\right)^l
= \exp\left(\frac{-1}{2c}\right) = C,\]
so we have to show that the factor involving $n$ makes no difference to
the limit. Now this factor is always less than~$1$, so the series
is absolutely convergent (and uniformly in $n$); so we can choose $r$ large
enough that the sum of $r$ terms of each sequence is close to its limit.
Then, since the factors tend to $1$ as $n\to\infty$, for $n$ large each of
these $r$ terms is close to its limit. So the assertion is true: that is,
$F'_1(n)\sim Cg(n)$.

\subparagraph{Case $F'_2$:} A permutation which has $k=n-l$ cycles and is
not an involution has at least $n-2l+1$ fixed points, and there are at most
\[{n\choose 2l-1}(2l-1)!=n(n-1)\cdots(n-2l+2)\]
such permutations. So, ignoring signs,
\begin{eqnarray*}
\frac{F'_2(n)}{g(n)}
&\le& \sum_{l\ge0}\frac{(n-l)(n-l-1)\cdots(n-2l+2)}{n(n-1)\cdots(n-l+1)}\,
c^{-l}\\
&<& \frac{1}{n}\cdot\frac{1}{1-c^{-1}},
\end{eqnarray*}
which is $O(1/n)$.

\subparagraph{Case $F'_3$:} We simply observe that there are at most $n!$
such permutations, so
\[\frac{F'_3(n)}{g(n)} \le n!\sum_{k=1}^{\lfloor n/2\rfloor}
\left(\frac{k!}{n!}\right)^2.\]
Now $n!/(k!)^2\ge{n\choose\lfloor n/2\rfloor}\ge(2-\epsilon)^n$ for large~$n$,
so this sum is $O((2-\epsilon)^{-n}n/2)=o(1)$ as $n\to\infty$.
\qed

\noindent {\bf Third Proof:}
If one is interested
in asymptotic enumeration of $F(n)$, (\ref{mobius}), being a double
sum over terms of alternating sign, is on first sight rather unsuitable
for an asymptotic analysis. The expression in Proposition \ref{prop2} is 
also an alternating sum.
We present a derivation of the asymptotic form
of $F(n)$ based on the following elegant and elementary identity, which
gives $F(n)$ as a sum of positive terms.
\begin{prop}
\begin{equation}
\label{ident1}
F(n)=\sum_{k=0}^\infty\sum_{l=0}^\infty\frac1{2^{k+l+2}}{kl\choose n}\;.
\end{equation}
\end{prop}
\proof
Insert
\begin{equation}\label{sum1}
1=\sum_{k=i}^\infty\frac1{2^{k+1}}{k\choose i}=
\sum_{l=j}^\infty\frac1{2^{l+1}}{l\choose j}
\end{equation} into (\ref{fsum}) and resum using (\ref{firstdisp}).\qed

We start the asymptotic analysis by rewriting (\ref{ident1}) as
\begin{equation}\label{centralsum}
n!F(n)=\sum_{k=0}^\infty\sum_{l=0}^\infty
\frac{k^n}{2^{k+1}}\frac{l^n}{2^{l+1}}\frac{(kl)_n}{(kl)^n}\;,
\end{equation}
where $(x)_n=x(x-1)\ldots(x-n+1)$ is the falling factorial. Given the identity
\begin{equation}\label{pseries}
P(n)=\sum_{k=0}^\infty\frac{k^n}{2^{k+1}}\;,
\end{equation}
which follows from expanding (\ref{pgenfunc}), (\ref{centralsum}) is bounded
above by $P(n)^2$, as the factor $(kl)_n/(kl)^n$ takes values in $[0,1]$. 

For $n\leq kl$, a straightforward expansion of the factor gives
\begin{equation}\nonumber
\frac{(kl)_n}{(kl)^n}=\exp\left[
-\sum_{j=1}^\infty\frac{B_{j+1}(n)-B_{j+1}(0)}{j(j+1)(kl)^j}\right]\;.
\end{equation}
Here, we have used that $\sum_{k=0}^{n-1}k^j=(B_{j+1}(n)-B_{j+1}(0))/(j+1)$ where $B_j(x)$
is a Bernoulli polynomial. It follows that
\begin{equation}\label{quotientasy}
\frac{(kl)_n}{(kl)^n}=e^{-\frac{n^2}{2kl}}\left(1+O(n/kl)+O(n^3/(kl)^2)\right)\;.
\end{equation}
(This argument will be presented more thoroughly for $(z)_n$ with complex-valued $z$ 
in the next section.)
The sum (\ref{centralsum}) is dominated by terms around $k=l=n/\log2$, so that we
expect the correction to give $e^{-(\log2)^2/2}$, which in turn would imply
$n!F(n)\sim P(n)^2e^{-(\log2)^2/2}$. The difference is given by
\begin{equation}\nonumber
n!F(n)-P(n)^2e^{-(\log2)^2/2}=\sum_{k=0}^\infty\sum_{l=0}^\infty
\frac{k^n}{2^{k+1}}\frac{l^n}{2^{l+1}}\left(\frac{(kl)_n}{(kl)^n}-e^{-(\log2)^2/2}\right)\;.
\end{equation}
To proceed we choose $m_0<n/\log2<m_1$ and 
split the summation. We obtain
\begin{eqnarray}\nonumber
\left|n!F(n)-P(n)^2e^{-(\log2)^2/2}\right|
&\leq&\sum_{k=m_0}^{m_1}\sum_{l=m_0}^{m_1}
\frac{k^n}{2^{k+1}}\frac{l^n}{2^{l+1}}\left|\frac{(kl)_n}{(kl)^n}-e^{-(\log2)^2/2}\right|\\
&&+2\sum_{k=0}^\infty\frac{k^n}{2^{k+1}}\left(\sum_{k=0}^{m_0-1}\label{splitsum}
\frac{k^n}{2^{k+1}}+\sum_{k=m_1+1}^\infty\frac{k^n}{2^{k+1}}\right)\;.
\end{eqnarray}
Specifying $m_0=n/\log2-cn^\delta$ and $m_1=n/\log2+cn^\delta$ for $1/2<\delta<1$ and $c>0$, we use (\ref{quotientasy}) to estimate
\begin{eqnarray*}
\frac{(kl)_n}{(kl)^n}-e^{-(\log2)^2/2}
&=&e^{-\frac{n^2}{2kl}}\left(1+O(n/kl)+O(n^3/(kl)^2)\right)-e^{-(\log2)^2/2}\\
&=&e^{-(\log2)^2/2}\left(1+O(n^{\delta-1})\right)\left(1+O(n^{-1})\right)-e^{-(\log2)^2/2}\\
&=&O(n^{\delta-1})
\end{eqnarray*}
for $m_0\leq k,l\leq m_1$. This allows us to bound the first term in (\ref{splitsum}) by
$P(n)^2O(n^{\delta-1})$. To get a bound on the second term, we utilize the following 
\begin{lem}\label{lemma}
\begin{itemize}
\item[(a)]
For $K,n\in\mathbb N$ and $K<n/\log2$, 
\begin{equation}
\sum_{k=0}^K\frac{k^n}{2^{k+1}}\leq\frac{K^n}{2^K}\frac{e^{n/K}}2\frac1{e^{n/K}-2}\;.
\end{equation}
\item[(b)]
For $K,n\in\mathbb N$ and $K>n/\log2$, 
\begin{equation}
\sum_{k=K+1}^\infty\frac{k^n}{2^{k+1}}\leq\frac{K^n}{2^K}\frac{e^{n/K}}2\frac1{2-e^{n/K}}\;.
\end{equation}
\end{itemize}
\end{lem}
\begin{pf}
Part (a) follows from the estimate
\begin{eqnarray*}
\sum_{k=0}^K\frac{k^n}{2^{k+1}}
&=&\frac{K^n}{2^{K+1}}\sum_{k=0}^K\frac{(k/K)^n}{2^{k-K}}
\leq\frac{K^n}{2^{K+1}}\sum_{k=0}^K2^k(1-k/K)^n
\leq\frac{K^n}{2^{K+1}}\sum_{k=0}^\infty\left(2e^{-n/K}\right)^k
\end{eqnarray*}
and part (b) similarly from
\begin{eqnarray*}
\sum_{k=K+1}^\infty\frac{k^n}{2^{k+1}}
&=&\frac{K^n}{2^{K+1}}\sum_{k=K+1}^\infty\frac{(k/K)^n}{2^{k-K}}
\leq\frac{K^n}{2^{K+1}}\sum_{k=1}^\infty\frac{(1+k/K)^n}{2^k}
\leq\frac{K^n}{2^{K+1}}\sum_{k=1}^\infty\left(\frac{e^{n/K}}2\right)^k\;.
\end{eqnarray*}
\qed
\end{pf}
For $K=n/\log2\mp cn^\delta$, we find
$$\frac{K^n}{2^K}\frac{e^{n/K}}2\frac1{|e^{n/K}-2|}=\frac{n^ne^{-n}}{(\log2)^n}
e^{-\alpha n^{2\delta-1}}O(n^{1-\delta})=P(n)O\left(e^{-\alpha n^{2\delta-1}}\right)$$
where $\alpha=c^2(\log2)^2/2$. Using Lemma \ref{lemma}, we therefore bound the second term in (\ref{splitsum}) by $P(n)^2O\left(\exp(-\alpha n^{2\delta-1})\right)$. Altogether we find
$$
n!F(n)-P(n)^2e^{-(\log2)^2/2}
=P(n)^2\left(O(n^{\delta-1})+O\left(\exp(-\alpha n^{2\delta-1})\right)\right)
$$
and as $1/2<\delta<1$, we have
$$
\lim_{n\rightarrow\infty}\frac{n!F(n)}{P(n)^2}=e^{-(\log2)^2/2}
$$
which completes the proof.
\qed

\section{The asymptotics of $m_{kl}(n)$}

In this section we present results on the number of incidence matrices with
specified numbers of rows and columns. To obtain the desired asymptotic
form of $m_{kl}(n)$ from eqn.\ (\ref{mobius}), we need to deal with the
challenge that summing over large terms with alternating signs can lead to
enormous cancellations. Fortunately, there is a standard trick using the
calculus of residues.
\begin{prop}
\begin{equation}
\label{ident2}
m_{kl}(n)=\frac{k!l!}{n!}\mathrm{Res}(
\frac{(st)_n}{(s)_{k+1}(t)_{l+1}};\;s=\infty,\;t=\infty)\;.
\end{equation}
\end{prop}
\begin{pf}
Using the fact that
$$\mathrm{Res}(\Gamma(s),\;s=-m)=\frac{(-1)^m}{m!},$$
we write
\begin{eqnarray*}
m_{kl}(n)&=&\frac{k!l!}{n!}(-1)^{k+l}\sum_{i=0}^k\sum_{j=0}^l
\frac{(-1)^i}{i!}\frac{(-1)^j}{j!}\frac{(ij)_n}{(k-i)!(l-j)!}\\
&=&\frac{k!l!}{n!}\frac{(-1)^{k+l}}{(2\pi i)^2}
\int_{{\cal C}_{[-k,0]}}ds\int_{{\cal C}_{[-l,0]}}dt\;
\frac{(st)_n}{(s+k)_{k+1}(t+l)_{l+1}}\\
&=&\frac{k!l!}{n!}\frac1{(2\pi i)^2}
\int_{{\cal C}_{[0,k]}}ds\int_{{\cal C}_{[0,l]}}dt\;
\frac{(st)_n}{(s)_{k+1}(t)_{l+1}}\;.
\end{eqnarray*}
Here, the contours ${\cal C}_{[a,b]}$ encircle the (real) interval $[a,b]$ counterclockwise. As the integrand is a rational function, the contour integrals can be
expressed as residues at infinity. \qed
\end{pf}
This formulation allows us to do the asymptotic analysis via saddle point
analysis of a contour integral. We consider the scaling behaviour of $m_{k,l}(n)$
as $n\rightarrow\infty$ with $k=\kappa n$ and $l=\lambda n$ for fixed $\lambda$, 
$\kappa$. As a preparation, we state the following Lemma.
\begin{lem}\label{fallingident}
Let $n\in\mathbb N$ and $z\in\mathbb C$ with $|z|>n$. Then
\begin{equation}\label{fallingident1}
(z)_n=z^n\exp\left[
-\sum_{j=1}^\infty\frac{B_{j+1}(n)-B_{j+1}(0)}{j(j+1)z^j}\right]\;.
\end{equation}
Moreover, we have the asymptotic expansion
\begin{eqnarray}\label{fallingident2}\nonumber
\log(z)_n&\sim&(z+1/2)\log z-(z-n+1/2)\log(z-n)-n\\
&&+\sum_{k=1}^\infty\frac{B_{2k}}{2k(2k-1)}\left(\frac1{z^{2k-1}}-\frac1{(z-n)^{2k-1}}\right)\;.
\end{eqnarray}
as $|z-n|$ and $|z|$ tend to infinity.
Here $B_n(x)$ is the $n$-th Bernoulli polynomial and $B_n=B_n(0)$ the $n$-th Bernoulli number.
\end{lem}
\proof
We write
$$(z)_n=z^n\prod_{k=0}^\infty\left(1-\frac kz\right)\;.$$ 
For $|z|>n$ we take logarithms and expand $\log(1-k/z)$ in $k/z$. Exchanging the order
of summation and using that
$$\sum_{k=0}^{n-1}k^j=(B_{j+1}(n)-B_{j+1}(0))/(j+1)$$
gives (\ref{fallingident1}). One can obtain (\ref{fallingident2}) by substituting 
$B_n(x)=\sum_{k=0}^n{n\choose k}B_kx^{n-k}$ and exchanging the order of summation again. As
the double sum here is not absolutely convergent, the resulting series (\ref{fallingident2})
cannot be expected to be convergent. Instead of labouring to prove that one still arrives at an 
asymptotic expansion, we point out that the result is just the difference between the 
Stirling series for $\log(z!)$ and $\log((z-n)!)$ and argue that standard arguments used in 
the derivation of the Stirling series also lead to (\ref{fallingident2}). In contrast 
with the Stirling series, which is valid for $|\arg(z)|<\pi$, the validity of this expansion is not restricted to a sector of the 
complex plane.
\qed

From (\ref{fallingident1}) we obtain that
\begin{equation}
\label{fallingasy1}
(z)_n=z^ne^{-\frac{n^2}{2z}}\left(1+O(n/|z|)+O(n^3/|z|^2)\right)\;,
\end{equation}
whereas from (\ref{fallingident2}) we obtain that
\begin{equation}
\label{fallingasy2}
(z)_n=(z-n)^{n-z-1/2}z^{z+1/2}e^{-n}\left(1+O(1/|z|)+O(1/|z-n|)\right)\;.
\end{equation}

We now state the main theorem of this section.
\begin{theorem}\label{main}
For fixed $\sigma,\tau>0$,
\begin{equation}
m_{kl}(n)\sim\frac{n^{2n}}{n!}
e^{nw(\sigma)}v(\sigma)e^{nw(\tau)}v(\tau)e^{-\frac1{2\sigma\tau}}
\end{equation}
with
\begin{eqnarray*}
w(x)&=&x(1-e^{-1/x})\log(1-e^{-1/x})+\log x-e^{-1/x}\;,\\
v(x)&=&\sqrt{\frac{x(1-e^{-1/x})}{x(1-e^{-1/x})-e^{-1/x}}}
\end{eqnarray*}
and $k=n\sigma(1-e^{-1/\sigma})$, $l=n\tau(1-e^{-1/\tau})$.
\end{theorem}

\begin{pf}
In order to evaluate the integral 
$$
\frac{n!}{k!l!}m_{kl}(n)=
\frac1{(2\pi i)^2}\int_{{\cal C}_{[0,k]}}ds\int_{{\cal C}_{[0,l]}}dt\;
\frac{(st)_n}{(s)_{k+1}(t)_{l+1}}
$$
asymptotically, we approximate the integrand uniformly using (\ref{fallingasy1}) and
(\ref{fallingasy2}) on contours satisfying $|s|=R_s>k$ and $|t|=R_t>l$. We find
$$
\frac{n!}{k!l!}m_{kl}(n)=
\frac{e^{k+l}}{(2\pi i)^2}\int_{|s|=R_s}ds\int_{|t|=R_t}dt\;
\frac{(s-k)^{s-k-1/2}}{s^{s-n+1/2}}
\frac{(t-l)^{t-l-1/2}}{t^{t-n+1/2}}
e^{-\frac{n^2}{2st}}(1+R)
$$
where
$$
R=O(n/R_sR_t)+O(n^3/(R_sR_t)^2)+O(1/R_s)+O(1/R_t)+O(1/|R_s-k|)+O(1/|R_t-l|)\;.
$$
Substituting $s=n\sigma$, $t=n\tau$, $k=n\kappa$, $l=n\lambda$, this simplifies to
$R=O(1/n)$ and we arrive at
$$
m_{kl}(n)\sim
\kappa^{n\kappa}\lambda^{n\lambda}n^ne^n\frac{\sqrt{2\pi\kappa\lambda n}}{(2\pi i)^2}\int_{|\sigma|=\rho_\sigma}d\sigma\int_{|\tau|=\rho_\tau}d\tau\;
e^{nf(\sigma,\kappa)}g(\sigma,\kappa)e^{nf(\tau,\lambda)}g(\tau,\lambda)e^{-\frac1{2\sigma\tau}}
$$
with $R_s=n\rho_\sigma$ and $R_t=n\rho_\tau$ and
$$f(x,y)=(x-y)\log(x-y)-(x-1)\log x\quad\mbox{and}\quad g(x,y)=((x-y)x)^{-1/2}\;.$$
As $n$ tends to infinity, each integration is dominated by saddle points 
$\sigma_s$ and $\tau_s$ on the positive real axis. The saddle point equation is
$$
0=\partial_1f(\sigma_s,\kappa)
=\log\left(1-\frac\kappa{\sigma_s}\right)+\frac1{\sigma_s}\;,
$$
(here, $\partial_1$ denotes taking the derivative with respect to the first argument)
and an identical expression for $\tau_s$.
There exists a unique positive solution $\sigma_s(\kappa)$, and a standard saddle-point evaluation (see e.g.
\cite{bleistein}) gives
$$
m_{kl}(n)\sim
\kappa^{n\kappa}\lambda^{n\lambda}n^ne^n\frac{\sqrt{2\pi\kappa\lambda n}}{(2\pi)^2}
e^{nf(\sigma_s,\kappa)}g(\sigma_s,\kappa)e^{nf(\tau_s,\lambda)}g(\tau_s,\lambda)e^{-\frac1{2\sigma_s\tau_s}}
\frac{2\pi}{\sqrt{\partial_1^2f(\sigma_s,\kappa)\partial_1^2f(\tau_s,\lambda)}}
$$
which simplifies to the desired result.
\qed
\end{pf}

Theorem \ref{main} can be used for a fourth proof of Theorem \ref{F1111}
via an asymptotic evaluation of the sum over $m_{kl}(n)$. The sum is
dominated by terms near $\sigma_s=\tau_s=1/\log2$ from whence it follows
that the distribution has a peak about $k_s=l_s=n/(2\log2)$.

We conclude this paper with giving an identity for $m_{kl}(n)$ 
which is a refinement of Proposition \ref{prop2}.
\begin{prop}\label{propprop}
\begin{equation}
\label{ident3}
m_{kl}(n)=\frac{k!l!}{n!}\sum_{r=1}^ns(n,r)S(r,k)S(r,l)
\end{equation}
\end{prop}
\proof

Using $(ij)_n=\sum_{r=1}^ns(n,r)(ij)^r$, we write (\ref{mobius}) as
$$
m_{kl}(n)=\frac1{n!}\sum_{r=1}^ns(n,r)\sum_{i\le k}\sum_{j\le l}(-1)^{k+l-i-j}{k\choose i}{l\choose j}
(ij)^r
$$
and resum using $k!S(r,k)=\sum_{i\le k}(-1)^{k-i}{k\choose i}i^r$.
\qed

Inversion of (\ref{ident3}) gives
\begin{equation}\label{inverseident3}
k!S(n,k)l!S(n,l)=\sum_{r=1}^nr!S(n,r)m_{kl}(r)
\end{equation}
which has a straightforward combinatorial interpretation, as $r!S(n,r)$ is the number
of preorders of $n$ elements into $k$ blocks. The left hand side of (\ref{inverseident3})
is just the number of ways of choosing two preorders of an $n$-set into $k$ and $l$ blocks,
respectively. The right hand side of (\ref{inverseident3})
corresponds to counting the number of ways in which elements of an $n$-set can be distributed 
into $r$ cells of a $k\times l$-array, where the cells are given by $k\times l$-incidence 
matrices with $r$ ones, for arbitrary $r$.

Summing (\ref{ident3}) over $k$ and $l$ provides another proof of Proposition
\ref{prop2}. Proposition \ref{propprop} could also be used as a basis for Theorem \ref{main}.
We leave this as an exercise for the reader.

\end{document}